\documentclass[10pt,draft,reqno]{amsart}
     \makeatletter
     \def\section{\@startsection{section}{1}%
     \z@{.7\linespacing\@plus\linespacing}{.5\linespacing}%
     {\bfseries
     \centering
     }}
     \def\@secnumfont{\bfseries}
     \makeatother
\usepackage{color}
\usepackage{graphicx}
\usepackage{subfigure}
\usepackage{epsfig}
\usepackage{cite}
\usepackage{amsmath, amssymb}
\usepackage{fancybox}
\usepackage{listings}
\usepackage{verbatim}
\usepackage{url}
\usepackage{esint}
\usepackage{fancyhdr}
\usepackage{booktabs}
\usepackage{caption}
\usepackage{dsfont}
\newtheorem{theorem}{Theorem}[section]
\newtheorem{lemma}[theorem]{Lemma}

\newtheorem{corollary}[theorem]{Corollary}

\usepackage{mathrsfs,amsthm,amsmath}

\usepackage{pstricks}
\usepackage{amssymb}
\usepackage{pgfplots}
\pgfplotsset{width=10cm,compat=1.9}

\usepackage{graphicx}
\usepackage[mathscr]{euscript}

\usepackage{bbm,color,bm}
\usepackage{hyperref}
 \usepackage{mathptmx}

\newcommand{\rmd}{{\mathrm d}}
\newcommand{\rme}{{\mathrm e}}
\newcommand{\la}{{\lambda}}

\newcommand{\ep}{{\varepsilon}}
\newcommand{\1}{\mathbbm{1}}
\newcommand{\pd}{{\partial}}

\newcommand{\sT}{{\mathcal T}}
\newcommand{\cP}{{\mathcal P}}
\newcommand{\sL}{{\mathcal L}}

\newcommand{\ttT}{{\mathtt T}}

\newcommand{\PP}{{\mathbb P}}

\newcommand{\RR}{{\mathbb R}}

\newcommand{\EE}{{\mathbb E}}
\newcommand{\TT}{{\mathbb T}}

\newcommand{\bp}{{\mathbf p}}

\begin{document}
\title{Metastability for telegraph processes in a double-well potential}
\footnote{Forthcoming in J.Appl.Probab.}
\author{ Nikita
Ratanov
}\thanks{Chelyabinsk State University,\;
rtnnkt@gmail.com}
\date{}
\maketitle
\begin{abstract}{
In this paper, we study equations with nonlinearity 
in the form of a double-well potential, randomised by a velocity-switching (telegraph) 
stochastic process.
If the speed parameters of the randomisation
are small, then this dynamics has one metastable uncertainty interval and two invariant attractors. 
The probabilities of leaving the metastable interval through the upper boundary are determined,
as well as characteristics of the first crossing times. 
Invariant measures are also found.

When and if the direction of the telegraph process velocity
 coincides with the direction of the periodic change in potential,
the system can go into a metastable state, having 
received a time window for the interwell transition.

The obtained results can be used as an alternative to stochastic resonance models.
}

\noindent\emph{Keywords:} piecewise deterministic process; first passage time; invariant measure; 
stochastic resonance\\
\end{abstract}
\section{Introduction}
Anomalous relativistic (\cite{Gzyl}, see also \cite{zigzag2}) 
diffusion processes with finite propagation speed 
often seems more natural for modelling than the classical diffusion approach. 
There are several different terminological and ideological systems for dealing 
with such classes of anomalous diffusions.
In physics, people prefer to say ``continuous time random walks" \cite{Weiss}, 
``flip-flop" \cite{Bremaud}, 
or ``zig-zag processes"  \cite{zigzag,zigzag2}.
In mathematics, the terms ``piecewise-deterministic processes", \cite{Davis}, 
or ``telegraph processes", \cite{Goldstein,Kac}, are commonly used. 
The study of classical stochastic differential equations 
with regime switching (and their applications)
is also widely represented, see for example \cite{BJPS,Mao}.

In this paper, we study a piecewise deterministic nonlinear system whose dynamics 
is specified by a double-well potential and a two-state Markov random process. It turns out that such 
a system can have metastable states as an important feature. 

This problem setting continues the author's previous studies on stochastic differential equations 
based on piecewise deterministic processes. These studies have recently begun with  
Ornstein-Uhlenbeck-like processes, where the telegraph process replaces
the Wiener process in the Langevin equation.
For such stochastic dynamics with its various modifications, invariant measures and
distributions of first passage times were obtained, 
see \cite{MathBio,MCAP2021,MCAP2022,GSA22,GSA23}. 
Here we develop similar ideas in a \emph{nonlinear} setting.

Dynamic systems subject to random perturbations are actively studied 
 in both physical and mathematical literature. 
In the simplest case, such systems are usually based on deterministic dynamics, 
which follows the ordinary differential equation 
$\rmd\gamma(t)=-U'(\gamma(t))\rmd t,$ $t>0,$ 
equipped with a random perturbation, where the function $U$ is the potential.

Let $x$ be a stationary point, i.e., $U'(x)=0.$ If $U''(x)>0$, then the position $x$ is asymptotically 
 stable, that is $\gamma(t)\to x$ as $t\to+\infty$ for any starting point near $x$.
On the contrary, if $U''(x)<0,$ then $x$ unstable and repulsive, i.e., any 
path that starts near $x$ moves away from $x$. 
\emph{The interplay between different kinds of stationary points 
under different randomisations of this equation
is of special interest.}
The first and most popular randomisation method is to add
a stochastic term, such as white noise, to the deterministic differential equation, 
see  \cite{imkeller-book}.
Later, other randomisation options have emerged, 
for example,  the L\'evy process \cite{y}, ``coloured" noise (e.g., the Ornstein-Uhlenbeck process) \cite{c,Wilkinson},
the fractional Gaussian \cite{goychuk}.

When the phases coexist, the phenomenon of 
hysteresis  may occur,  see the discussion in \cite{das,y,verma}.
An introduction to metastability can be found in \cite{Bovier}.

The simplest case of one-dimensional setting has recently been generalised in
a number of works on metastability effects arising in
the case of diffusion processes that are stopped or reflecting on the boundary of domain,
see \cite{Koralov-AP}, and more generally, on metastability
 at and near the surfaces that remain invariant under diffusion processes, \cite{Koralov}. 

When using various randomisation methods, stochastic resonance can occur, 
where a small change in the input signal causes a large output signal in the system.
Such dynamics, with a transition from one potential well to another, 
can be interpreted to qualitatively explain glacial cycles in the history of the Earth 
(from a warm state to an ice age and back), see
\cite{imkeller-book}. Similar phenomena occur in various applied settings,
from biology to metallurgy, \cite{Strogatz}.
For example, this approach is effective when modelling a single neurone, 
\cite{bulsara,tuckwell,nicolis,MathBio} (see also a recent review of these models in 
\cite{Carrillo}), and artificial neuronal sets
\cite{simonito}, as well as dynamics of competing populations \cite{Wissel} and
chemical reaction systems \cite{Vellela}. 
See also detailed overviews of other applications in \cite{rmp,McD,intech}.

The idea of stochastic resonance is successfully used for example by adding white noise 
$\sigma(T)\rmd W(t)$ to a deterministic differential equation,   
where $W$ is standard Brownian motion and the coefficient $\sigma(T)$
is chosen to \emph{maximise the coefficient of spectral power amplification}. 
See the seminal works  \cite{benzi,nicolis}.
The second approach involves the so-called 'effective dynamics',
 based on a two-state Markov chain with a time-periodic infinitesimal generator,
which artificially ensures a spatial interwell transition,  see e.g. \cite{imkeller-book}. 

Our randomisation methodology consists in replacing 
the diffusion term $W(t)$ by a relativistic (finite velocity) process $\TT(t)$.
The idea of such replacement becomes more and more popular in physics, \cite{Gzyl,Sandev}.
With this randomisation, the model receives certain time windows for interwell transitions 
instead of stochastic resonance options. 
In this settings, we present an infinitesimal generator of the process (Theorem \ref{theo:generator}), 
compare the probabilities of which of the attractor the process 
will end up on when the velocities are small
(Theorem \ref{theo:uncertain}), and derive formulae for the average first passage times 
when the process has already left the metastable position
(Theorem \ref{theo:startingNearAttractor}). 
Invariant measures inside the attractors are also found (Theorem \ref{theo:sm}).

The paper is organised as follows. The next section is devoted to the problem setting 
and description of the structure of the underlying dynamics. In Section \ref{sec3} we analyse 
the distributions of passage times. Section \ref{sec4} presents invariant distributions.
\section{Preliminaries and problem statement}
Let $\TT(t)$ be a telegraph process with alternating velocities $c_0, c_1,\;c_0>c_1$,
\begin{equation}
\label{def:T}
\TT(t)=\int_0^tc_{\xi(s)}\rmd s,\qquad t>0,
\end{equation}
where $\xi=\xi(t)\in\{0, 1\}$ is a two-state Markov process with  
alternating switching intensities $\la_0, \;\la_1$, \cite{RK}.
Consider a stochastic process $X(t)=X^x(t)$ defined by the equation
\begin{equation}\label{eq:X}
\mathrm d X(t)=-U'(X(t))\rmd t+\rmd \TT(t),\qquad t>0,
\end{equation}
and the initial condition $X(0)=x$. 
Process $X=X(t)$ can be considered as a telegraph process $\TT$ in the potential $U$.

Our previous studies \cite{MathBio,MCAP2021,MCAP2022,GSA22,GSA23} were based on 
the so-called Kac-Ornstein-Uhlenbeck processes,
which are defined by the equation
\begin{equation*}
\rmd \bar X(t)=\left(c_{\xi(t)}-\gamma_{\xi(t)}\bar X(t)\right)\rmd t,\qquad t>0,
\end{equation*}
i.e. the case of alternating single-well potentials 
$U_i=\gamma_ix^2/2-c_ix$.

For the quasilinear setting \eqref{eq:X}, suppose that the function $U(x),\;U\in C^2(-\infty, \infty),$ is a \emph{ double-well potential}, i.e.
it has a local maximum at the point $x_0$, two local minima at the points 
$x_{-},\;x_+,\;x_{-}<x_0<x_{+},$ and $U''(x)\ne0,\;x\in\{x_-, x_0, x_+\}.$  
As a typical example, one can use the symmetric function
$U(x)=\dfrac{x^4}{4}-\dfrac{x^2}{2}$.

The process $X(t)$ is continuous,
piecewise continuously differentiable and is formed by the sequential alternation of two 
deterministic patterns $\gamma_0(t, x)$ and $\gamma_1(t, x),$ which follow the equations
\begin{equation}
\label{eq:gamma}\aligned
\frac{\pd\gamma_0(t, x)}{\pd t}=c_0-U'(\gamma_0(t, x)),&\qquad
\frac{\pd\gamma_1(t, x)}{\pd t}=c_1-U'(\gamma_1(t, x)),\qquad t>0,\\
&\gamma_0(0, x)=\gamma_1(0, x)=x.
\endaligned
\end{equation}
One pattern is replaced by another 
at random times $\tau_n,$ when the underlying Markov process $\xi(t)$ switches.

The initial value problems \eqref{eq:gamma} are equivalent to the pair of integral equations
\begin{equation}
\label{eq:patterns-int}
\int_x^{\gamma_0(t, x)}\frac{\rmd z}{c_0-U'(z)}=t,\qquad \text{and}\qquad
\int_x^{\gamma_1(t, x)}\frac{\rmd z}{c_1-U'(z)}=t,\qquad t\ge0,
\end{equation}
%
which can be rewritten in the equivalent form:
\begin{equation}
\label{eq:gammaPhi}
\Phi_0(\gamma_0(t, x))=\Phi_0(x)+t\qquad \text{and}\qquad
\Phi_1(\gamma_1(t, x))=\Phi_1(x)+t.
\end{equation}
Here 
 $\Phi_0$ and $\Phi_1$ are rectifying diffeomorphisms  
defined up to an arbitrary additive constant, by the equations
\begin{equation}
\label{def:Phii}
\Phi_0'(y)=\frac{1}{c_0-U'(y)},\qquad \Phi_1'(y)=\frac{1}{c_1-U'(y)}.
\end{equation}
Notice that the function $\Phi_i=\Phi_i(y)$ is continuous and monotone 
between the critical points of the potential $U_i=U(x)-c_ix$.

To avoid blow-up, assume that the integral $\int_{x}\frac{\rmd y}{U_i'(y)}$ diverges 
at any critical point of $U_i:$
\begin{equation}
\label{eq:int-infty}
\int_{x}\frac{\rmd y}{U_i'(y)}=\infty,\qquad x\in\{x_-, x_0, x_+\}.
\end{equation}

From \eqref{eq:patterns-int} it follows  that the pattern $\gamma_i(t, x)$
arising from the state  $(x, i),\;i\in\{0, 1\},$ 
also satisfies the initial value problem
\begin{equation}
\label{eq:Lgamma}\aligned
\frac{\pd\gamma_i}{\pd t}(t, x)=&L_i^x[\gamma_i(t, x)],\qquad t>0,\\
\gamma_i(t, x)|_{t\downarrow0}=&x,
\endaligned\end{equation}
where $L_i^x=-U_i'(x)\dfrac{\pd}{\pd x}=(c_i-U'(x))\dfrac{\pd}{\pd x},\; i\in\{0, 1\}.$

The properties of the process $X$ depend significantly on the velocity amplitudes $c_0$ and $c_1$.

 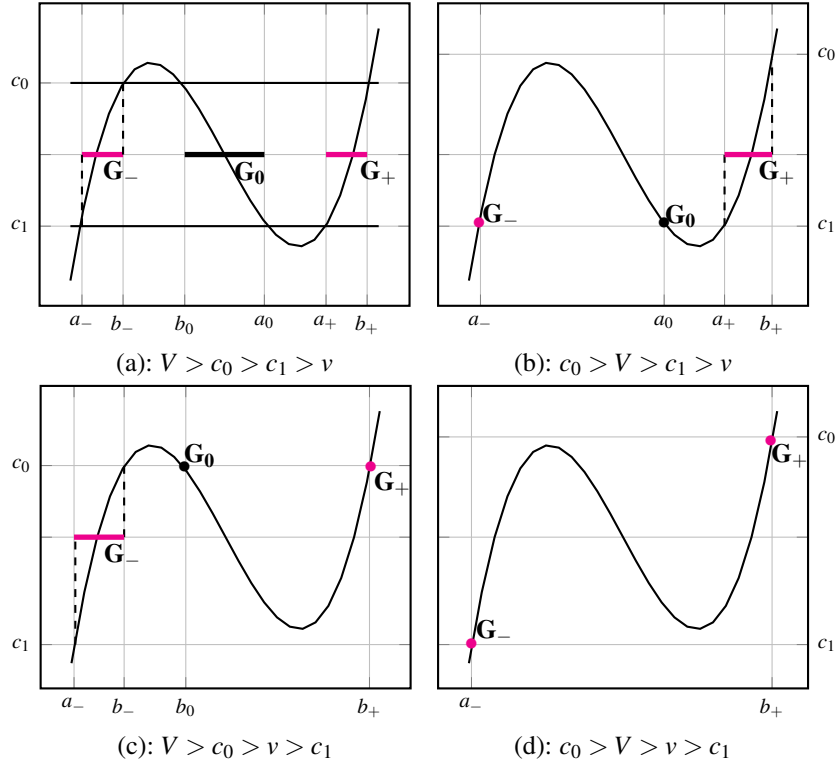
\begin{figure}[t]
\pgfplotsset{
small
}
\begin{center}
\begin{tabular}{rl}
\begin{tikzpicture}[baseline,trim axis left,domain=-1.2:1.2]
\begin{axis}[thick,grid=both,
no markers,
xtick={-1.11,-0.79,-0.31,0.31,0.79,1.11},
ytick={-0.3,0,0.3},
yticklabels={$c_1$,,$c_0$},
xticklabels={$a_-$,$b_-$,$b_0$,$a_0$,$a_+$,$b_+$},
xlabel={(a):\;$V>c_0>c_1>v$},xlabel style={font=\normalsize},
xmajorgrids,
 ymajorgrids,
]
\addplot [draw=black]{x^3-x};
\addplot[magenta,line width=2pt] coordinates {(-1.11,0) (-0.79,0)};
\addplot[line width=2pt] coordinates {(-0.31,0) (0.31,0)};
\addplot[magenta,line width=2pt] coordinates {(0.79,0) (1.11,0)};
\addplot[thick] coordinates {(-1.2, 0.3) (1.2, 0.3)};
\addplot[thick] coordinates {(-1.2, -0.3) (1.2, -0.3)};
\addplot[dashed] coordinates {(-1.11, -0.3) (-1.11, 0)};
\addplot[dashed] coordinates {(-0.79, 0) (-0.79, 0.3)};
      \draw  node at (40,450) {$\mathbf{G_-}$};
          \draw  node at (240,450) {$\mathbf{G_+}$};
		 \draw  node at (140,450) {$\mathbf{G_0}$};
\end{axis}
\end{tikzpicture}
&
\begin{tikzpicture}[baseline,trim axis right,domain=-1.2:1.2]
\begin{axis}[thick,no markers,
yticklabel pos=upper,
xtick={-1.11,0.32,0.79,1.16},
ytick={-0.3,0,0.42},
yticklabels={$c_1$,,$c_0$},
xticklabels={$a_-$,$a_0$,$a_+$,$b_+$},
xlabel={(b):\;$c_0>V>c_1>v$},xlabel style={font=\normalsize},
xmajorgrids,
 ymajorgrids
]
\addplot [draw=black,thick] {x^3-x};

\addplot[dashed] coordinates {(0.79, -0.3) (0.79, 0)};
\addplot[dashed] coordinates {(1.16, 0) (1.16, 0.4)};
\addplot[magenta,line width=2pt] coordinates {(0.79,0) (1.16,0)};
 \draw  node at (240,450) {$\mathbf{G_+}$};
 \draw node at (152,240) {$\bullet$}; 
  \draw node at (165,270) {$\mathbf{G_0}$}; 
   \draw[magenta] node at (8,240) {$\bullet$};
  \draw node at (25,270) {$\mathbf{G_-}$}; 
\end{axis}
\end{tikzpicture}
\\
\begin{tikzpicture}[baseline,trim axis right,domain=-1.2:1.2]
\begin{axis}[thick,
no markers,
xtick={-1.18,-0.79,-0.31,1.12},
ytick={-0.45,0,0.3},
yticklabels={$c_1$,,$c_0$},
xticklabels={$a_-$,$b_-$,$b_0$,$b_+$},
xlabel={(c):\;$V>c_0>v>c_1$},xlabel style={font=\normalsize},
xmajorgrids,
 ymajorgrids,
]
\addplot [draw=black]{x^3-x};
\addplot[dashed] coordinates {(-1.17, -0.45) (-1.17, 0)};
\addplot[dashed] coordinates {(-0.79, 0) (-0.79, 0.3)};
\addplot[magenta,line width=2pt] coordinates {(-1.18,0) (-0.79,0)};
 \draw node at (88,820) {$\bullet$}; 
  \draw node at (100,870) {$\mathbf{G_0}$}; 
  \draw node at (40,450) {$\mathbf{G_-}$}; 
   \draw[magenta] node at (233,820) {$\bullet$};
          \draw  node at (250,750) {$\mathbf{G_+}$};
\end{axis}

\end{tikzpicture}
&
\begin{tikzpicture}[baseline,trim axis right,domain=-1.2:1.2]
\begin{axis}[thick,no markers,
yticklabel pos=upper,
xtick={-1.18,1.16},
ytick={-0.45,0,0.42},
yticklabels={$c_1$,,$c_0$},
xticklabels={$a_-$,$b_+$},
xlabel={(d):\;$c_0>V>v>c_1$},xlabel style={font=\normalsize},
xmajorgrids,
 ymajorgrids
]
\addplot [draw=black,thick] {x^3-x};

 \draw[magenta] node at (235,932) {$\bullet$}; 
  \draw node at (250,870) {$\mathbf{G_+}$}; 
   \draw[magenta] node at (2,80) {$\bullet$};
  \draw node at (20,130) {$\mathbf{G_-}$}; 
\end{axis}
\end{tikzpicture}
\\
\end{tabular}
\end{center}
   \caption{$U'(x)$ and $G_0,\;G_\pm$.}
    \label{figU'}
\end{figure}

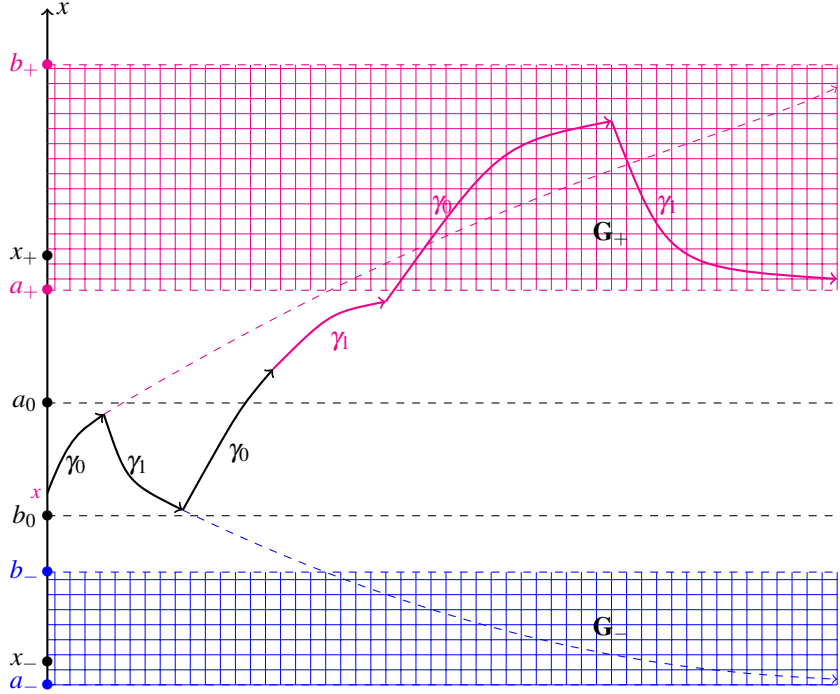
\begin{figure}[th]
\begin{tikzpicture}[x=1.5cm,y=1.5cm,domain=0:8]  
\draw[thick,->] (-5,-2) -- (-5,4)                 node[right] {$x$};                
     \draw[magenta
     ]  node at (-5.1,-0.3) {\footnotesize$x$};
       \draw (-5,1.8) node {$\bullet$} node at (-5.2,1.8) {$x_+$};
         \draw (-5,-1.8) node {$\bullet$} node at (-5.2,-1.8) {$x_-$};
         \draw node at (-4.75,-0.05) {$\gamma_0$};
         \draw node at (-4.2,-0.05) {$\gamma_1$};
          \draw node at (-3.3,0.05) {$\gamma_0$};
           \draw[magenta] node at (-2.4,1.05) {$\gamma_1$};
           \draw[magenta] node at (-1.5,2.25) {$\gamma_0$};
           \draw[magenta] node at (0.5,2.25) {$\gamma_1$};
    \draw (-5,0.5) node {$\bullet$} node at (-5.2,0.5) {$a_0$};
  \draw (-5,-0.5) node {$\bullet$} node at (-5.2,-0.5) {$b_0$};
   \draw[magenta
   ] (-5,1.5) node {$\bullet$} node at (-5.2,1.5) {$a_+$};
   \draw[magenta
   ] (-5,3.5) node {$\bullet$} node at (-5.2,3.5) {$b_+$};
     \draw[blue
     ] (-5,-1) node {$\bullet$} node at (-5.2,-1) {$b_-$};
   \draw[blue
   ] (-5,-2) node {$\bullet$} node at (-5.2,-2) {$a_-$};
\draw[magenta,
dashed] (-5,1.5) -- (2,1.5);   
\draw[magenta,
dashed] (-5,3.5) -- (2,3.5);  
\draw[blue,
dashed] (-5,-1) -- (2,-1);   
\draw[blue,
dashed] (-5,-2) -- (2,-2);  
\draw[dashed] (-5,0.5) -- (2,0.5);  
\draw[dashed] (-5,-0.5) -- (2,-0.5);   
\draw[step=0.2cm,magenta,
very thin ] (-5,3.5) grid (2,1.5);
\draw[step=0.2cm,blue,
very thin ] (-5,-1) grid (2,-2);
\draw node at (0,-1.5) {$\mathbf {G_-}$};
\draw node at (0,2) {$\mathbf {G_+}$};
\draw[thick, ->] (-5,-0.3) .. controls (-4.8,0.23) and (-4.7,0.25)  .. (-4.5,0.4);
\draw[thin,dashed,magenta,
->] (-4.5,0.4)  .. controls (-1,2.5) and (1,2.8)  .. (2,3.3);
\draw[thick,->] (-4.5,0.4) .. controls (-4.3,-0.2) ..(-3.8,-0.45);
\draw[thin,dashed,blue,
->] (-3.8,-0.45)  .. controls (-1,-1.8) and (1,-1.9)  .. (2,-1.95);
\draw[thick,->] (-3.8,-0.45).. controls (-3.3, 0.45)  .. (-3.0, 0.8);
\draw[thick,magenta,
->] (-3.0, 0.8).. controls(-2.5,1.3)..(-2,1.4);
\draw[thick,magenta,
->] (-2,1.4).. controls (-1, 2.8).. (0,3);
\draw[thick,magenta,
->] (0, 3).. controls (0.5,1.7).. (2,1.6);
            \end{tikzpicture}
             \caption{Sample path of $X$, case A.}
       \label{figfig}
\end{figure}  

\textbf{A: An unstable uncertainty interval and two attractors. Fig.\ref{figU'}(a).}

Let both parameters $c_0$ and $c_1$ be small enough so that $V>c_0>c_1>v$,
where $V$ and $v$ are, respectively, a local maximum and a local minimum of $U'$. 
In this case, both
$U_0(x)=U(x)-c_0x$ and $U_1(x)=U(x)-c_1x$ are still double-well potentials with alternation between them.
Let $a_0, a_\pm$ be the critical points of $U_1,\;a_-<a_0<a_+,$
and $b_0, b_\pm$ be 
 the critical points of $U_0,\;,\;b_-<b_0<b_+,$ see Fig.\ref{figU'}(a).
 
Note that $X(t)\equiv a_-,\;X(t)\equiv a_+$ and $X(t)\equiv a_0$ are  
stationary solutions of \eqref{eq:X} under the initial state $\xi(0)=1$. Pattern $\gamma_1=\gamma_1(t, x)$
is repelled by $a_0$, while being attracted to $a_-$ or $a_+$.
Respectively, $X(t)\equiv b_-,\;X(t)\equiv b_+$ and $X(t)\equiv b_0$ are %
stationary solutions under $\xi(0)=0,$ and
pattern $\gamma_0$
is repelled by $b_0,$ being attracted to points $b_\pm$.

As a result, the process $X$ has two attractors: 
$G_-=[a_-,\; b_-]\subset(-\infty, b_0),\;
G_+=[a_+, \;b_+]\subset(a_0, +\infty)$ and a metastable set $G_0=(b_0, a_0)$.
The process $X^x,$ starting from the metastable set $G_0,\;x\in G_0,$  can temporarily oscillate inside, 
 but once it leaves this interval through the upper or lower bound,  it never returns, 
and then forever remains above $a_0$ or below $b_0,$ respectively, being attracted to $G_+$ or $G_-$.
See the sample of a path in Fig.\ref{figfig}.

If the process $X=X^x(t)$ starts from $x=X^x(0)>a_0$, then after a finite transition time it
ends up in $G_+$, and when starting from $x<b_0$ the process for a.s. finite time falls into $G_-.$

Note that both sets $G_-$ and  $G_+$ are invariant under dynamics:
\begin{equation}
\label{eq:inv}
\PP\{X(t)\in G_-\;\forall  t>0~|~X(0)\in G_-\}=1,\qquad 
\PP\{X(t)\in G_+\;\forall  t>0~|~X(0)\in G_+\}=1.
\end{equation}

Figure \ref{figPhipm} shows plots of $\Phi_0$ and $\Phi_1$ in this case.

\textbf{B: Alternating single-well and double-well potentials, only one attractor. Fig.\ref{figU'}(b), (c).}

When $c_0$ becomes large, so that $c_0>V$, the potential $U_0(x)$ 
has only one critical point $b_+$ becoming single-well,
  while $U_1$ is still double-well, 
and has three critical points. In this case, the process $X$ has only one attractor $G_+$,
see Fig.\ref{figU'}(b). After some transition period, with probability 1 $X$ ends up in $G_+$. 
A symmetric situation occurs when $c_1$ becomes less than $v$, see Fig.\ref{figU'}(c).

\textbf{C: Alternating single-well potentials, one attractor. Fig.\ref{figU'}(d).}

In the case $c_0>V>v>c_1,$ both potentials $U_0$ and $U_1$ become single-well. 
The dynamics has only one attractor $G=(a_-, b_+).$
When the process $X$ begins in $x=X(0)\in G$, it randomly fluctuates between the 
attractive levels $a_-$ and $b_+$. If it starts outside, then it ends at $G$ after some transition.

Two other less interesting single-well situations with a single attractor, are not shown in Fig.\ref{figU'}. 
They correspond to $c_0>c_1>V>v$ and 
$c_1<c_0<v<V$. See Fig.\ref{figc0c1} for a diagram of the different regimes.

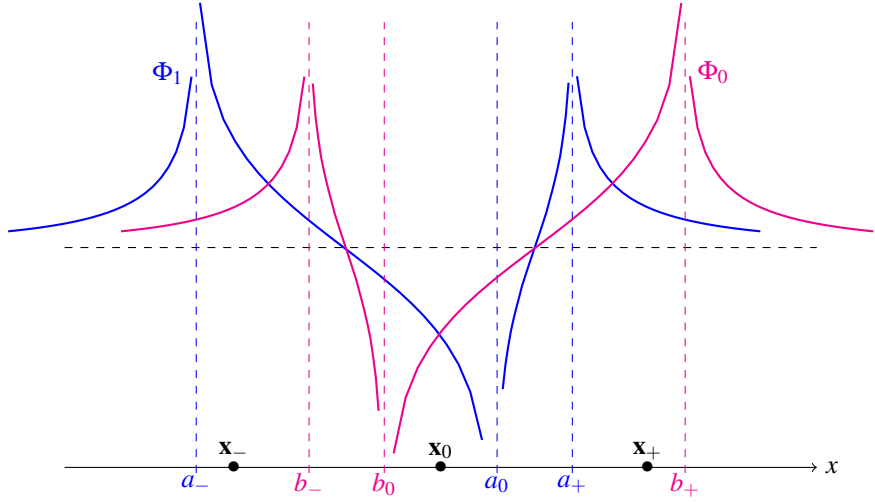
\begin{figure}[t]
\begin{tikzpicture}[x=2.5cm,y=1.5cm,domain=-4.7:4.7]  
\draw[dashed] (-2,0) -- (2,0)
          node at (1.1,-1.8) {$\mathbf x_+$}
      node at (-1.1,-1.8) {$\mathbf x_-$}
      node at (0, -1.8) {$\mathbf x_0$};      
                \draw (1.1,-1.95) node {$\bullet$};
            \draw (-1.1,-1.95) node {$\bullet$};
            \draw (0,-1.95)  node {$\bullet$};    
    \draw[blue,
    dashed] (0.3,-2) -- (0.3,2)   
     node at (0.3, -2.1) {$a_0$};
    \draw[blue,
    thick,domain=-2.3:-1.325]
    plot(\x,{-0.5*ln(1-(\x+0.3)^(-2))});
     \draw[blue,
     thick,domain=-1.28:0.22]
    plot(\x,{-0.5*(ln(1+(\x+0.3))-ln(0.57-(\x+0.3)))})
     node at (-1.45,1.55) { $\Phi_1$};
       \draw[blue,
                   dashed] (-1.3,-2) -- (-1.3,2)
     node at (-1.3, -2.1) {$a_-$};
  \draw[blue,
  thick,domain=0.33:0.68]
  plot(\x,{-0.5*(ln(-(\x+0.3)+1)-ln(\x-0.3)});
     \draw[blue,
     thick,domain=0.725:1.7]
    plot(\x,{-0.5*ln(1-(\x+0.3)^(-2))});
     \draw[blue,
     dashed] (0.7,-2) -- (0.7,2)
     node at (0.7, -2.1) {$a_+$};
    \draw[->] (-2,-1.95) -- (2,-1.95)
                 node[right] {$x$};           
             \draw[magenta,
             dashed] (-0.3,-2) -- (-0.3,2)
    node at (-0.3, -2.1) {$b_0$};
    \draw[magenta,
    thick,domain=-1.7:-0.725]
    plot(\x,{-0.5*ln(1-(\x-0.3)^(-2))});
          \draw[magenta,
          thick,domain=-0.68:-0.33]
    plot(\x,{-0.5*(ln((\x+1.7)-1)-ln(-0.31-\x)}) 
    node at (1.45,1.55) {\magenta
    $\Phi_0$};
          \draw[magenta,
          thick,domain=-0.25:1.28]
    plot(\x,{-0.5*(ln(-(\x-0.3)+1)-ln(0.29+\x))});
     \draw[magenta,
     dashed] (-0.7,-2) -- (-0.7,2)
     node at (-0.7, -2.1) {$b_-$};
     \draw[magenta,
     thick,domain=1.325:2.3]
    plot(\x,{-0.5*ln(1-(\x-0.3)^(-2))});
     \draw[magenta,
     dashed] (1.3,-2) -- (1.3,2)
     node at (1.3, -2.1) {$b_+$};
            \end{tikzpicture}
        \caption
    {Functions $\Phi_0(x)$ and $\Phi_1(x),$ case A.
      }
    \label{figPhipm}
\end{figure}
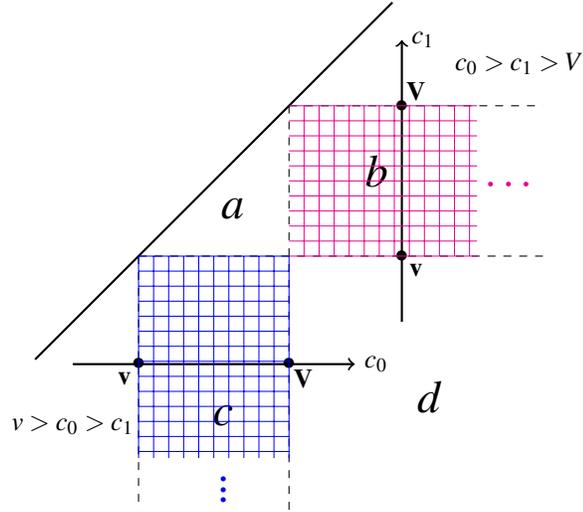
\begin{figure}[th]
\begin{center}
\begin{tikzpicture}[x=1.25cm,y=1.25cm,domain=-4.7:4.7] 
    \draw[thick,->] (-1.5,-1.95) -- (1.5,-1.95)       node[right] {$c_0$};   
 \draw[thick,->] (2,-1.5) -- (2,1.5)       node[right] {$c_1$};   
\draw[thick] (-1.9,-1.9) -- (1.9,1.9);
     \draw (-0.8,-1.95) node {$\bullet$};
       \draw (0.8,-1.95) node {$\bullet$};
      \draw (2,-0.8) node {$\bullet$};
           \draw (2,0.8) node {$\bullet$};
           \draw[dashed] (-0.8,-0.8) --(3.5,-0.8) node at (-0.95, -2.1) {$\mathbf v$};
      \draw[dashed] (-0.8,-0.8) -- (-0.8,-3.5) node at (2.15, -0.95) {$\mathbf v$};
      \draw[dashed] (0.8,-3.5) -- (0.8,0.8) node at (0.95, -2.1) {$\mathbf V$};
      \draw[dashed] (0.8,0.8) -- (3.5,0.8) node at (2.15, 0.95) {$\mathbf V$};
      \draw node at (3.25,1.25) {$c_0>c_1>V$};
            \draw node at (-1.5,-2.6) {$v>c_0>c_1$};
             \draw node at (0.2,-0.3) {{\huge$a$}};
              \draw node at (2.5,0.1) {{\huge$b\qquad\magenta\ldots$}};
               \draw node at (2.3,-2.3) {{\huge$d$}};
               \draw node at (0.1,-2.5) {{\huge$c$}};
               \draw node at (0.1,-3.2) {\blue\huge$\vdots$};
               \draw[step=0.2cm,magenta,very thin ] (0.8,-0.8) grid (2.8,0.8);
               \draw[step=0.2cm,blue,very thin ] (-0.8,-2.95) grid (0.8,-0.8);
           \end{tikzpicture}
        \caption
    {Scheme of different modes. Triangle $\mathbf{a},\;V>c_0>c_1>v$ (case A);
 half-bands $\mathbf{b}$ and $\mathbf{c},$ correspond to $c_0>V>c_1>v$ 
 and $V>c_0>v>c_1$ respectively (case B);
quadrant $\mathbf{d},\;c_0>V>v>c_1$ (case C)
 (designations according to Fig.\ref{figU'} ).   }
    \label{figc0c1}
    \end{center}
\end{figure}
The process $\Xi=(X^x(t),\;\xi(t))$ is a Markov process with the state space $\RR\times\{0, 1\}.$

Let $\tau^0$ be the time of first switching from state 0 to state1,
$\tau^1$ corresponds to the reverse switching.
Due to the Markov property, the following equalities in distribution hold:
\begin{equation}
\label{eq:conditioning}
\begin{aligned}
  \big[X^x(t)~|~\xi(0)=0\big]  & \stackrel{D}{=} \gamma_0(t, x)\cdot\1_{\{\tau^0>t\}}
  +\1_{\{\tau^0<t\}}\cdot\big[X^{\gamma_0(\tau^0, x)}(t-\tau^0)~|~\xi(0)=1\big], \\
    \big[X^x(t)~|~\xi(0)=1\big]  & \stackrel{D}{=} \gamma_1(t, x)\cdot\1_{\{\tau^1>t\}}
  +\1_{\{\tau^1<t\}}\cdot\big[X^{\gamma_1(\tau^1, x)}(t-\tau^1)~|~\xi(0)=0\big].
\end{aligned}
\end{equation}
Based on \eqref{eq:conditioning} one can obtain an infinitesimal generator of $\Xi$.
\begin{theorem}\label{theo:generator}
The infinitesimal generator of $\Xi$ is given by $\Lambda+\sL,$ where 
\begin{equation}
\label{def:L}
\sL=\sL^x
=\begin{pmatrix}
  L_0^x    &   0 \\
  0    &  L_1^x
\end{pmatrix},
\end{equation}
and $
\Lambda=\begin{pmatrix}
   -\la_0 & \la_0   \\
   \la_1   &  -\la_1
\end{pmatrix}$ is the generator of the two-state Markov process $\xi.$
\end{theorem}
See the proof in Appendix.
\section{Passage times}\label{sec3}
Let $\sT(x, y)$ be the time when the process $X(t)$ 
first reaches a given threshold $y$, starting from point $x$,
\[
\sT(x, y)=\inf\{t>0~|~X(0)=x,\;X(t)=y\}.
\]
The distribution of this random variable 
is determined by the moment generating function 
$\bm u=(u_0(q, x, y),\; u_1(q, x, y))^\ttT,$  
\begin{equation*}
\aligned
u_0(q, x, y)=&\EE[\exp(-q\sT(x, y))~|~\xi(0)=0],\\ 
u_1(q, x, y)=&\EE[\exp(-q\sT(x, y))~|~\xi(0)=1],
\endaligned\qquad q>0.\end{equation*}

By $t_0^*(x, y)$ and $t_1^*(x, y)$ we denote  
the time to reach level $y$ without switching states, 
starting from $(x, 0)$ and $(x, 1)$, respectively.
The value $t_i^*(x, y)$ is finite if 
 the equation
  $\gamma_i(t, x)=y$ is satisfied for some positive $t,\;t=t_i^*(x, y)$.
In this case we say that level $y$ is reachable from state $(x, i)$.  
Otherwise we set $t_i^*(x, y)=\infty$.
    
  If $y$ is reachable from the state $(x, i)$, then both
  $x$ and $y$ belong to the domain of a common continuous piece of $\Phi_i,$ and, moreover,
  since \eqref{eq:gammaPhi} holds,
  $\Phi_i(y)>\Phi_i(x)$ and $t_i^*(x, y)=\Phi_i(y)-\Phi_i(x).$ 
Moreover, $\lim\limits_{x\to y}t_i^*(x, y)=0.$

By virtue of
\eqref{def:Phii},   
\begin{equation}
\label{eq:Lt*}
L_i^xt_i^*(x, y)=-L_i^x\Phi_i(x)=-1,\qquad \text{when $t_i^*(x, y)$ is finite},\quad i\in\{0, 1\}.
\end{equation}

Let $y$ be a level reachable from both $(x, 0)$ and $(x, 1)$.
Conditioning on the first switching, like in \eqref{eq:conditioning}, 
we obtain a system of integral equations:
\begin{equation}
\label{eq:u-int}
\left\{
\begin{aligned}
   u_0(q, x, y)=\rme^{-(q+\la_0)t_0^*(x, y)}+ & 
   \int_0^{t_0^*(x, y)}\la_0\rme^{-(q+\la_0)\tau}u_1(q, \gamma_0
   (\tau, x), y)\rmd\tau,\\
   u_1(q, x, y)=\rme^{-(q+\la_1)t_1^*(x, y)}+&
   \int_0^{t_1^*(x, y)}\la_1\rme^{-(q+\la_1)\tau}u_0(q, \gamma_1
   (\tau, x), y)\rmd\tau. 
\end{aligned}
\right.
\end{equation}
On the contrary, in the case $t_i^*(x, y)=+\infty,$ 
 the corresponding integral equation of \eqref{eq:u-int} can be interpreted as
\begin{equation}
\label{eq:u-int-infty}
u_i(q, x, y)=\int_0^\infty\la_i\rme^{-(q+\la_i)\tau}u_{1-i}(q, \gamma_i(\tau, x), y)\rmd\tau,\qquad i\in \{0, 1\}.
\end{equation}

\begin{theorem}\label{theo:LL}
The moment generating function $\bm u=(u_0(q, x, y),\;u_1(q, x, y))^\ttT$ of 
the first passage time $\sT(x, y)$
satisfies the equation\textup:
\begin{equation}
\label{eq:u-diff-matrix}
\big(\Lambda+\sL)\bm u(q, x, y)=q\bm u(q, x, y),\qquad q>0,
\end{equation}
where $\Lambda+\sL$ is the generator of Markov process $\Xi.$ 

Both entries $u_0(q, x, y)$ and $u_1(q, x, y)$ of $\bm u$ separately
obey the telegraph  \textup{(\cite{borodin})} equations of the form\textup,  
\begin{equation}
\label{eq:equ0equ1}\aligned
L_1^xL_0^xu_0(q, x, y)-\big[(q+\la_1)L_0^x+(q+\la_0)L_1^x\big]u_0(q, x, y)
+q(q+\la_0+\la_1)u_0(q, x, y)&=0,\\
L_0^xL_1^xu_1(q, x, y)-\big[(q+\la_1)L_0^x+(q+\la_0)L_1^x\big]u_1(q, x, y)
+q(q+\la_0+\la_1)u_1(q, x, y)&=0.
\endaligned\end{equation}
\end{theorem}
Equations  \eqref{eq:u-diff-matrix} and \eqref{eq:equ0equ1} are provided with boundary conditions
that depend on whether the level $y$ is reachable from $(x, 0)$ and $(x, 1)$.
\proof
Apply  $L_0^x$ and $L_1^x$ to the first and the second equation 
of the system \eqref{eq:u-int}, respectively.

First, let the threshold $y$ be reachable from $(x, 0)$. 
By virtue of \eqref{eq:Lgamma} and \eqref{eq:Lt*},
the first equation of \eqref{eq:u-int} yields 
\[\aligned
L_0^xu_0(q, x, y)=(q+\la_0)\rme^{-(q+\la_0)t_0^*(x, y)}
-&\la_0\rme^{-(q+\la_0)t_0^*(x, y)}u_1(q, \gamma_0(t_0^*(x, y), x), y)\\
+&\int_0^{t_0^*(x, y)}\la_0\rme^{-(q+\la_0)\tau}\frac{\pd}{\pd\tau}
\left[u_1(q, \gamma_0(\tau, x), y)\right]\rmd\tau.
\endaligned\]
Integrating by parts, one can see that
\[\aligned
L_0^xu_0(q, x, y)=(q+\la_0)\rme^{-(q+\la_0)t_0^*(x, y)}-&\la_0u_1(q, x, y)\\
+&(q+\la_0)\int_0^{t_0^*(x, y)}\la_0\rme^{-(q+\la_0)\tau}u_1(q, \gamma_0(\tau, x), y)\rmd\tau.
\endaligned\]

Proceeding similarly with the second equation of \eqref{eq:u-int},
both in the case of finite and infinite $t_1^*(x, y)$, 
as a result, we obtain the system
\begin{equation}
\label{eq:u-diff}
\left\{
\begin{aligned}
L_0^xu_0(q, x, y)=    & (q+\la_0)u_0(q, x, y)-\la_0u_1(q, x, y), \\
 L_1^xu_1(q, x, y)=   &  -\la_1u_0(q, x, y)+(q+\la_1)u_1(q, x, y),
\end{aligned}
\right.\end{equation}
which  in matrix form coincides with equation \eqref{eq:u-diff-matrix}.

Equations \eqref{eq:equ0equ1} follow from \eqref{eq:u-diff} by eliminating variables.  
Indeed, applying $L_1^x$ to 
the first equation of system \eqref{eq:u-diff}, we get
$
L_1^xL_0^xu_0(q, x, y)=(q+\la_0)L_1^xu_0(q, x, y)-\la_0L_1^xu_1(q, x, y).
$ 
From the second equation of \eqref{eq:u-diff} one can obtain 
\[
L_1^xL_0^xu_0(q, x, y)
=(q+\la_0)L_1^xu_0(q, x, y)+\la_0\la_1u_0(q, x, y)-\la_0(q+\la_1)u_1(q, x, y).
\]  
The first equation of  \eqref{eq:u-diff}  can be rewritten as
$\la_0u_1(q, x, y)=(q+\la_0)u_0(q, x, y)-L_0^xu_0(q, x, y),$ which gives the first equation of
 \eqref{eq:equ0equ1}; the second equation is obtained similarly. 
\endproof
If $\PP\{\sT(x, y)<\infty\}=1$, then 
the mean value  
$\bm m(x, y)=(m_0(x, y),\; m_1(x, y))^\ttT,$ 
\[ m_0(x, y)=\EE[\sT(x,y)~|~\xi(0)=0], \qquad m_1(x, y)=\EE[\sT(x,y)~|~\xi(0)=1],\]
can be found by differentiation,
\begin{equation}
\label{eq:m0m1=}
\left.m_0(x, y)=-\frac{\pd u_0(q, x, y)}{\pd q}\right|_{q=0},\qquad 
\left.m_1(x, y)=-\frac{\pd u_1(q, x, y)}{\pd q}\right|_{q=0}.
\end{equation}

Theorem \ref{theo:LL} implies the following result, general for Markov processes.
\begin{corollary}
Let the mean value of $\sT(x, y)$ be finite\textup.
%
Then 
 the mean vector $\bm m(x, y)$ satisfies the equation
  \begin{equation}
\label{eq:m0m1-1-diff}
(\Lambda+\sL^x)\bm m(x, y)=-\bm1=(1, 1)^\ttT.
\end{equation}
\end{corollary}

Let us first consider the case of small $c_0$ and $c_1$,
i.e. $V>c_0>c_1>v$,
see Fig.\ref{figU'}(a).

\subsection{Small $c_0$ and $c_1$: probabilities of leaving the uncertainty  interval}
In the case $V>c_0>c_1>v$,  
let the initial point $x$ be inside the interval
$G_0$, i.e. $X(0)=x\in G_0=(b_0, a_0),$  case A, Fig.\ref{figU'}(a).
 
 Since both $a_0$ and $b_0$ are repulsive,
the trajectory that has once crossed one of these thresholds
  never returns  to the interval $(b_0, a_0)$.
We are interested in the conditional probabilities of exiting the interval 
  $G_0=(b_0, a_0)$ through the upper bound $a_0,$
\begin{equation}
\label{def:p0p1}
p_0(x)=\PP\{\sT(x, a_0)<\infty~|~\xi(0)=0\},\qquad 
p_1(x)=\PP\{\sT(x, a_0)<\infty~|~\xi(0)=1\}.
\end{equation}
The moment generating functions $u_0(q, x, a_0)$ and $u_1(q, x, a_0),\;x\in G_0,$ 
  are increasing. 
 Further, $0<u_0(q, x, a_0)<1,\;0<u_1(q, x, a_0)<1$ for $b_0<x<a_0.$
 Moreover, 
\[\aligned
\lim_{q\downarrow0}\left[u_0(q, x, a_0)\right]=&\PP\left\{\sT(x, a_0)<\infty\;\left.\right|\;\ep(0)=0\right\}=p_0(x)\\
&\qquad\text{and}\qquad \\
\lim_{q\downarrow0}\left[u_1(q, x, a_0)\right]=&\PP\left\{\sT(x, a_0)<\infty\;\left.\right|\;\ep(0)=1\right\}=p_1(x)\endaligned\]  
are the probabilities to leave the interval $(b_0, a_0)$ through the upper bound $a_0$, 
starting from the position
$X(0)=x$ with the initial state $\xi(0)=0$ and $\xi(0)=1,$ respectively.

The probabilities $p_0(x)$ and $p_1(x)$ 
depend on the following constants,
\begin{equation}
\label{def:B0B1} 
B_0=\int_{b_0}^{a_0}\frac{\Psi(y)}{c_0-U'(y)}\rmd y,\qquad 
B_1=\int_{b_0}^{a_0}\frac{\Psi(y)}{c_1-U'(y)}\rmd y.
\end{equation}
Here $\Psi(y)=\exp\left(\la_0\Phi_0(y)+\la_1\Phi_1(y)\right).$ 

The integrals in \eqref{def:B0B1} converge, since
\[
\left|\frac{\Psi(y)}{c_i-U'(y)}\right|\le const\cdot\exp\left(\la_i\Phi_i(y)\right)
\cdot\left|\Phi_i'(y)\right|,\qquad y\in G_0,
\]
and 
\[
\int_{b_0}^{a_0}\rme^{\la_i\Phi_i(y)}\left|\Phi_i'(y)\right|\rmd y
=\left|\exp\left(\la_i\Phi_i(y)\right)\bigg|_{y=b_0+}^{y=a_0-}\right|<\infty.
\]
See Fig.\ref{figPhipm}.

  \begin{theorem}\label{theo:uncertain}
Let $X(0)=x,\;x\in(b_0,\;a_0)$.
The probabilities $p_0(x)$ and $p_1(x)$ defined by \eqref{def:p0p1}
   have the form\textup:
  \begin{align}
\label{eq:uncertainp0}
p_0(x)=&\PP\{\sT(x, a_0)<\infty~|~\xi(0)=0\}
=1-B_0^{-1}\int_x^{a_0}\frac{\exp\left(\la_0\Phi_0(y)+\la_1\Phi_1(y)\right)}{c_0-U'(y)}\rmd y,
\\
p_1(x)=&\PP\{\sT(x, a_0)<\infty~|~\xi(0)=1\}
=B_1^{-1}\int_{b_0}^x\frac{\exp\left(\la_0\Phi_0(y)+\la_1\Phi_1(y)\right)}{c_1-U'(y)}\rmd y,
\label{eq:uncertainp1}
\end{align}
where $B_0, B_1$ are defined by \eqref{def:B0B1}.
  \end{theorem}
  \proof
By definition it follows that
 \begin{equation}
\label{eq:bc0-u0u1}
\left.p_0(x)\right|_{x\uparrow a_0}=1,\qquad
\left.p_1(x)\right|_{x\downarrow b_0}=0,
\qquad q>0.
\end{equation}

  Conditioning on the first switching like in \eqref{eq:conditioning}, similarly to \eqref{eq:u-int},
we obtain
\begin{equation}
\label{eq:u-intG0}
u_0(0, x, a_0)=p_0(x)=\rme^{-\la_0t_0^*(x, a_0)}
+\int_0^{t_0^*(x, a_0)}\la_0\rme^{-\la_0\tau}p_1(\gamma_0(\tau, x))\rmd\tau.
\end{equation}
Since the process $X$ leaving the interval $G_0$ through the lower bound
forever remains below $b_0$, then instead of \eqref{eq:u-intG0} and \eqref{eq:u-int-infty} we obtain
\begin{equation}
\label{eq:u-intG1}
u_1(0, x, a_0)=p_1(x)=\int_0^{t_1^*(x, b_0)}\la_1\rme^{-\la_1\tau}
p_0(\gamma_1(\tau, x))\rmd\tau.
\end{equation}
Here $t_0^*(x, a_0)=\Phi_0(a_0)-\Phi_0(x)$ and $t_1^*(x, b_0)=\Phi_1(b_0)-\Phi_1(x);$
function $\Phi_0$ is increasing, and function $\Phi_1$ is decreasing on $(b_0, a_0).$

From \eqref{eq:u-intG0}-\eqref{eq:u-intG1}
in the usual way by applying $L_0^x$ and $L_1^x$ followed by integration by parts,
differential equations similar to  \eqref{eq:equ0equ1} and \eqref{eq:u-diff}  are derived:
first \[(\Lambda+\sL)\bp(x)=\bm0,\qquad\bp(x)=(p_0(x), p_1(x))^\ttT,\] and then
\begin{equation}\label{eq:L0L1p}
\begin{aligned}
  L_1^xL_0^x p_0(x) & -(\la_1L_0^x+\la_0L_1^x)p_0(x)=0,  \\
  L_0^xL_1^x p_1(x) & -(\la_1L_0^x+\la_0L_1^x)p_1(x)=0.
\end{aligned}
\end{equation}
See Theorem \ref{theo:LL}.

Let 
\begin{equation}
\label{def:vv}
v_0(x)=L_0^x[p_0](x),\qquad v_1(x)=L_1^x[p_1](x).
\end{equation}
Note that the functions $p_0(x)$ and $p_1(x)$ increase, $x\in G_0.$

Since 
$c_1<U'(x)<c_0,$ 
then  $v_0(x)>0$ and $v_1(x)<0, \forall x\in G_0.$ 
In these notations, equations \eqref{eq:L0L1p} take a form
\[\aligned
(c_1-U'(x))\frac{\pd v_0}{\pd x}(x)-\la_1v_0(x)-\la_0\frac{c_1-U'(x)}{c_0-U'(x)}v_0(x)=&0,\\
(c_0-U'(x))\frac{\pd v_1}{\pd x}(x)-\la_1\frac{c_0-U'(x)}{c_1-U'(x)}v_1(x)-\la_0v_1(x)=&0,
\endaligned\]
which can be rewritten as
\begin{equation}
\label{eq:v0v1}\aligned
\frac{\rmd v_0}{\rmd x}(x)=\psi(x)v_0(x),&\qquad \frac{\rmd v_1}{\rmd x}(x)=\psi(x)v_1(x),\\
\qquad &b_0<x<a_0,
\endaligned\end{equation}
where 
\begin{equation}
\label{def:K}
\psi(x)
=\frac{\la_0}{c_0-U'(x)}+\frac{\la_1}{c_1-U'(x)}.
\end{equation}
From \eqref{eq:v0v1} and \eqref{def:K} by definition \eqref{def:Phii}  we obtain 
\begin{equation*}
v_0(x)=A_0\exp\left(\la_0\Phi_0(x)+\la_1\Phi_1(x)\right),\qquad 
v_1(x)=A_1\exp\left(\la_0\Phi_0(x)+\la_1\Phi_1(x)\right)
\end{equation*}
with indefinite constants $A_0$ and $A_1$.

Taking into account boundary conditions \eqref{eq:bc0-u0u1}, 
we obtain
\[
    p_0(x)=-\int_x^{a_0}\frac{v_0(y)}{c_0-U'(y)}\rmd y+1,   \qquad
    p_1(x)=\int_{b_0}^x  \frac{v_1(y)}{c_1-U'(y)}\rmd y.
\]

By the no-blow-up condition \eqref{eq:int-infty}, one can see that $A_0=B_0^{-1}$ and $A_1=B_1^{-1}$,
which gives the solution to equations \eqref{def:vv} in the form \eqref{eq:uncertainp0}-\eqref{eq:uncertainp1}.
\endproof
\subsection{Starting near attractor} 
Let the initial point $x$ and the threshold $y$ be located near  (or inside) the attractor. 
%

We study only the case of the upper half-line, $x, y>a_0$, 
and the attractor $G_+$, see Fig.\ref{figU'}(a) and (b). 
It is sufficient because the case of $x, y<b_0$ and the attractor $G_-$, 
Fig.\ref{figU'}(a) and (c), is symmetric.

Since both levels $a_+$ and $b_+$ are attractive, and $a_0,\; b_0$ are repulsive,  
 the process $X=X(t)$ increases a.s. as moving along the set $(a_0,\;a_+)$ 
 and decreases as moving along the set $[b_+,\;+\infty).$ Therefore,
 \[
 \left.\sT(x, y)\right|_{x>y\ge a_0}=\infty,\qquad
  \left.\sT(x, y)\right|_{b_+\le x<y}=\infty.\qquad a.s.
 \]
 
Inside the attractor $G_+$, the process $X=X(t)$ randomly oscillates between the attractive levels $a_+$ and $b_+$.
 
 The mean value  
$\bm m(x, y)=(m_0(x, y),\; m_1(x, y))^\ttT,$
can be expressed through the following integrals:
\begin{align}
\label{def:I}
  I_0(x, y)&=\int_x^y\frac{\beta(z_0, y)\rmd z_0}{c_0-U'(z_0)},     \qquad
  I_1(x, y)=\int_x^y\frac{\beta(z_1, y)\rmd z_1}{c_1-U'(z_1)},  \\
    \label{def:J0}
   J_0(x, y) & 
   = \int_{\Delta_0(x, y)}\frac{\beta(z_0, z_1)}{(c_0-U'(z_0))(c_1-U'(z_1))}\rmd z_0\rmd z_1,\\
       \label{def:J1}
   J_1(x, y) & 
   = \int_{\Delta_1(x, y)}\frac{\beta(z_1, z_0)}{(c_0-U'(z_0))(c_1-U'(z_1))}\rmd z_0\rmd z_1.
\end{align}
Here $\beta(x, y)=\dfrac{\Psi(x)}{\Psi(y)}
=\exp\left(\la_0(\Phi_0(x)-\Phi_0(y))+\la_1(\Phi_1(x)-\Phi_1(y))\right),$ 
and
$\Delta_0(x, y), \Delta_1(x, y)\subset\RR^2$ are two triangles defined by
\[\aligned\Delta_0(x, y)=&\{(z_0, z_1)~|~x<z_0<y,\;z_0<z_1<y\},\\
\Delta_1(x, y)=&\{(z_0, z_1)~|~x<z_1<y,\;z_1<z_0<y\}.\endaligned\]

According to the above definitions, 
 in both cases: $a_0<x<y\le a_+$ (with $U'(x), U'(y)<c_1<c_0$)
 and $x>y\ge b_+$  (with $U'(x), U'(y)>c_0>c_1$),
 the integrals $I_0(x, y), I_1(x, y)$ and $J_0(x, y), J_1(x, y)$ converge. 

\begin{theorem}\label{theo:startingNearAttractor}
\begin{itemize}
\item
Let $x, y>a_0,$  see Fig.\ref{figU'}(a) and (b).

The mean values $m_0(x, y)$ and $m_1(x, y)$ are given explicitly as follows.

\begin{itemize}
  \item In both cases\textup, $x<y<a_+$ and $x>y>b_+,$ the following formulae are valid\textup:
  \begin{align}
\label{eq:m0-sol-1}
    m_0(x, y)&
  =I_0(x, y)+(\la_0+\la_1)J_0(x, y),\\
    \label{eq:m1-sol-1}
     m_1(x, y)&
     =I_1(x, y)+(\la_0+\la_1)J_1(x, y).
         \end{align}
    
    \item If $x<y$ and $y\in[a_+, b_+),$ then  
  \begin{align}
\label{eq:m0-sol-4}
    m_0(x, y)&
=(1+\la_0/\la_1)I_0(x, y)+(\la_0+\la_1)J_0(x, y), \\
    \label{eq:m1-sol-4}
     m_1(x, y)&
   =\frac{1}{\la_1}+(\la_0+\la_1)J_1(x, y).
 \end{align}
  \item If $x>y$ and $y\in(a_+, b_+],$ then
  \begin{align}
\label{eq:m0-sol-5}
    m_0(x, y)&
    =\frac{1}{\la_0}+(\la_0+\la_1)J_0(x, b_+),
 \\
    \label{eq:m1-sol-6}
     m_1(x, y)&
   =(1+\la_0/\la_1)I_1(x, b_+)+(\la_0+\la_1)J_1(x, b_+).
 \end{align}
   \end{itemize}
Here $I_0,\;I_1$ and $J_0,\;J_1$
are defined by \eqref{def:I}-\eqref{def:J1}.
\item
Let $x, y<b_0,$  see Fig.\ref{figU'}(a) and (c).

The mean values $m_0(x, y)$ and $m_1(x, y)$ are determined by the following formulae\textup, 
symmetrical to those given above.

\begin{itemize}
  \item In both cases\textup, $x>y>b_-$ and $x<y<a_-,$ the following formulae are valid\textup:
  \begin{align}
\label{eq:m0-sol-2}
    m_0(x, y)&
  =I_0(x, y)+(\la_0+\la_1)J_0(x, y),\\
    \label{eq:m1-sol-2}
     m_1(x, y)&
     =I_1(x, y)+(\la_0+\la_1)J_1(x, y).
         \end{align}
    \item If $x>y$ and $y\in[a_-,b_-),$ then  
  \begin{align}
\label{eq:m0-sol-4-2}
     m_0(x, y)&
   =\frac{1}{\la_0}+(\la_0+\la_1)J_0(x, y),\\
    \label{eq:m1-sol-4-2}
    m_1(x, y)&
=(1+\la_1/\la_0)I_1(x, y)+(\la_0+\la_1)J_1(x, y);
 \end{align}
  \item If $x<y$ and $y\in(a_-,  b_-],$ then
  \begin{align}
\label{eq:m0-sol-5-2}
     m_0(x, y)&
   =(1+\la_1/\la_0)I_0(x, b_+)+(\la_0+\la_1)J_0(x, b_+),\\
    \label{eq:m1-sol-6-2}
    m_1(x, y)&
    =\frac{1}{\la_1}+(\la_0+\la_1)J_1(x, b_+).
 \end{align}
   \end{itemize}
   \end{itemize}
\end{theorem}
Notice that $X$ is structured in such a way that if
this process begins with $x\in(a_0,\;a_+]$ or with $x\in[b_+,\;+\infty),$
(symmetrically, $x\in[b_-, b_0)$ or $x\in(-\infty, a_-]$)
then with probability 1 it ends up inside the set $G_+$ ($G_-$)
 and remains there forever, randomly fluctuating between two attracting thresholds 
 $a_+$ and $b_+$.
  
In the rest of the paper, all proofs for the lower half-line ($x<b_0$) are symmetric to the proofs for $x>a_0,$
so we present only the latter.

To prove Theorem \ref{theo:startingNearAttractor},
we need some auxiliary statement, which is important in itself.
\begin{lemma}\label{theo:eqm}
Under the conditions of Theorem \ref{theo:startingNearAttractor}\textup,
the mean value of $\sT(x, y)$ is finite\textup, and the entries 
$m_0(x, y)$ and $m_1(x, y)$
of the mean vector $\bm m(x, y)$ are positive solutions of equation \eqref{eq:m0m1-1-diff}.

Furthermore\textup,
\begin{itemize}
\item
If $a_0<x<y<a_+$ $(b_0>x>y>b_-)$ or $x>y>b_+$ $(x<y<a_-),$ then \eqref{eq:m0m1-1-diff} is supplied with
zero boundary condition $\bm m(x, y)|_{x\to y}=\bm0.$
\item  If $x<y$ and $y\in[a_+, b_+),$ then 
\begin{equation}
\label{eq:bc1}
m_0(x, y)|_{x\to y}=0,\qquad m_1(x, y)|_{x\to y}=1/\la_1;
\end{equation}
if $x>y$ and $y\in(a_-, b_-],$ then
\begin{equation}
\label{eq:bc1-}
m_0(x, y)|_{x\to y}=1/\la_0,\qquad m_1(x, y)|_{x\to y}=0.
\end{equation}
\item If $x>y$ and $y\in(a_+, b_+],$  then
\begin{equation}
\label{eq:bc2}
m_0(x, y)|_{x\to y}=1/\la_0,\qquad m_1(x, y)|_{x\to y}=0;
\end{equation}
if $x<y$ and $y\in[a_-, b_-),$ then
\begin{equation}
\label{eq:bc2-}
m_0(x, y)|_{x\to y}=0,\qquad m_1(x, y)|_{x\to y}=1/\la_1.
\end{equation}
\end{itemize}

Moreover\textup, $m_0(x, y)$ and $m_1(x, y)$ satisfy separate second-order equations
\begin{equation}
\label{eq:m-diff3}
\left\{
\begin{aligned}
 L_1^xL_0^xm_0(x, y)- \big(\la_1L_0^x+\la_0L_1^x\big)m_0(x, y)=\la_0+\la_1,  &   \\
   L_0^xL_1^xm_1(x, y)- \big(\la_1L_0^x+\la_0L_1^x\big)m_1(x, y)=\la_0+\la_1.  &  
\end{aligned}
\right.
\end{equation}
\end{lemma}
\proof
First note that if $a_0<x<y<a_+$ or $x>y>b_+$, then
\emph{
  $t_0^*(x, y)$ and $t_1^*(x, y)$ are both finite},  
  $t_0^*(x, y)|_{x\to y}=t_1^*(x, y)|_{x\to y}=0,$ 
  functions $u_0,\;u_1$ follow the coupled integral equations equations \eqref{eq:u-int}, 
  $\sT(x, y)<\infty$ a.s.
  and $u_0(0, x, y)=u_1(0, x, y)=1.$ Differentiating in \eqref{eq:u-int} with respect to $q$ 
and then setting $q=0,$  we obtain the integral equations
 \begin{equation}
\label{eq:m-int1}
\left\{
\begin{aligned}
m_0(x, y)  =t_0^*(x, y)&\rme^{-\la_0t_0^*(x, y)}\\
+&\int_0^{t_0^*(x, y)}\la_0\tau\rme^{-\la_0\tau}\rmd\tau
+\int_0^{t_0^*(x, y)}\la_0\rme^{-\la_0\tau}m_1(\gamma_0(\tau, x), y)\rmd\tau,   \\
m_1(x, y)  =t_1^*(x, y)&\rme^{-\la_1t_1^*(x, y)}\\
+&\int_0^{t_1^*(x, y)}\la_1\tau\rme^{-\la_1\tau}\rmd\tau
+\int_0^{t_1^*(x, y)}\la_1\rme^{-\la_1\tau}m_0(\gamma_1(\tau, x), y)\rmd\tau.
\end{aligned}
\right.
\end{equation}
After simplification, system \eqref{eq:m-int1} becomes 
  \begin{equation}
\label{eq:m-int2}
\left\{
\begin{aligned}
m_0(x, y)  &=\frac{1-\rme^{-\la_0t_0^*(x, y)}}{\la_0} 
+\int_0^{t_0^*(x, y)}\la_0\rme^{-\la_0\tau}m_1(\gamma_0(\tau, x), y)\rmd\tau,   \\
m_1(x, y)  &=\frac{1-\rme^{-\la_1t_1^*(x, y)}}{\la_1}
+\int_0^{t_1^*(x, y)}\la_1\rme^{-\la_1\tau}m_0(\gamma_1(\tau, x), y)\rmd\tau.
\end{aligned}
\right.
\end{equation}


Note that in the cases A and B, 
the interval $G_+$ attracts both patterns, $\gamma_0(t, x)$ and $\gamma_1(t, x)$.
Therefore, 
 $t_i^*(x, y)|_{x\to y}=0$ if $x<y<a_+$  or $x>y>b_+$.
By \eqref{eq:m-int2}, this yields the zero boundary condition for \eqref{eq:m0m1-1-diff}:
$m_i(x, y)|_{x\to y}=0,\;i\in\{0, 1\}$. 

If $x<y=a_+$, then $t_0^*(x,y)$ is still finite, $\left.t_0^*(x, y)\right|_{x\uparrow y}=0$
and $t_1^*(x, y)=\infty$.
In this case, the first equations in \eqref{eq:m-int1} and \eqref{eq:m-int2} remain the same.
Hence, $\left.m_0(x, y)\right|_{x\uparrow y}=0$. 
The second equation in \eqref{eq:m-int2} takes the form
\begin{equation*}
m_1(x, y)=\frac{1}{\la_1}%
+\int_0^\infty\la_1\rme^{-\la_1\tau}m_0(\gamma_1(\tau, x), y)\rmd\tau,
\end{equation*}
which gives the same differential equation as before and the boundary condition \eqref{eq:bc1}.

If $x>y=b_+$, then similar reasoning again leads to the boundary condition 
\eqref{eq:bc2}.

By eliminating variables (as in Theorem \ref{theo:LL})
one can see  that
$m_0(x, y)$ and $m_1(x, y)$ satisfy two separate second-order equations
\eqref{eq:m-diff3}.
\endproof

\textit{The proof of Theorem \ref{theo:startingNearAttractor}} is similar to the proof of 
Theorem \ref{theo:uncertain}. It consists of two steps.
First, we transform the system \eqref{eq:m-diff3} of the second-order equations 
into a system of first-order equations for the functions
\begin{equation}\label{def:v0v1}
v_0(x, y):=L_0^xm_0(x, y),\qquad v_1(x, y):=L_1^xm_1(x, y).
\end{equation}
Note that \eqref{def:v0v1} can be solved in the form:
\begin{equation}
\label{eq:mv1}\aligned
m_0(x, y)=&\left.m_0(x, y)\right|_{x=y}-\int_x^y\frac{v_0(z_0, y)}{c_0-U'(z_0)}\rmd z_0,\\
m_1(x, y)=&\left.m_1(x, y)\right|_{x=y}-\int_x^y\frac{v_1(z_1, y)}{c_1-U'(z_1)}\rmd z_1,
\endaligned\end{equation}

Second, formulae \eqref{eq:mv1} will give the answer if we find $v_0$ and $v_1$. 

\emph{1st step}.
By definition of $L_0^x$ and $L_1^x$,
\[
L_0^xm_1(x, y)=\frac{c_0-U'(x)}{c_1-U'(x)}\cdot v_1(x, y),\qquad 
L_1^xm_0(x, y)=\frac{c_1-U'(x)}{c_0-U'(x)}\cdot v_0(x, y).
\]
Applying the change of variables \eqref{def:v0v1} to \eqref{eq:m-diff3}, we get 
\[\aligned
(c_1-U'(x))\frac{\pd v_0(x, y)}{\pd x}-\la_1v_0(x, y)-\la_0\frac{c_1-U'(x)}{c_0-U'(x)}v_0(x, y)
&=\la_0+\la_1,
\\
(c_0-U'(x))\frac{\pd v_1(x, y)}{\pd x}-\la_1\frac{c_0-U'(x)}{c_1-U'(x)}v_1(x, y)-\la_0v_1(x, y)
&=\la_0+\la_1.
\endaligned\]
Thus, system \eqref{eq:m-diff3} is equivalent to
\begin{equation}
\label{eq:v}
\left\{
\begin{aligned}
    \frac{\pd v_0}{\pd x}(x, y)& =\psi(x)v_0(x, y)+\frac{\la_0+\la_1}{c_1-U'(x)},  \\
   \frac{\pd v_1}{\pd x}(x, y)& =\psi(x)v_1(x, y)+\frac{\la_0+\la_1}{c_0-U'(x)},
\end{aligned}
\right.
\end{equation}
where $\psi(x)$ is defined by \eqref{def:K}, 
$\psi(x)=\dfrac{\la_0}{c_0-U'(x)}+\dfrac{\la_1}{c_1-U'(x)}=\Psi'(x)/\Psi(x)$.
The solution of \eqref{eq:v} is determined by the boundary conditions 
that follow from Lemma \ref{theo:eqm}.

\emph{2d step}.
Consider the case $a_0<x<y<a_+,$ with
$
   U'(x), U'(y) <c_1<c_0.    
 $ 
The case $x>y>b_+$ is symmetric with $U'(x),\; U'(y)>c_0>c_1.$ 

Since in this case, $\bm m(x, y)|_{x\to y}=\bm0$
(Lemma \ref{theo:eqm}),
then by virtue of \eqref{eq:m0m1-1-diff}, system \eqref{eq:v} is supplied with the boundary conditions
\begin{equation}
\label{eq:bcv}
v_0(x, y)|_{x\uparrow y}=v_1(x, y)|_{x\uparrow y}=-1.
\end{equation}
Therefore,
\begin{align}
\label{eq:v0-sol}
  v_0(z_0, y)  &\nonumber
  =\Psi(z_0)\Psi(y)^{-1}\bigg[-1-(\la_0+\la_1)\int_{z_0}^y\frac{\Psi(y)\Psi(z_1)^{-1}}{c_1-U'(z_1)}\rmd z_1\bigg]\\
 & =-\beta(z_0, y)-(\la_0+\la_1)\int_{z_0}^y\frac{\beta(z_0, z_1)}{c_1-U'(z_1)}\rmd z_1,\\
  \label{eq:v1-sol}
  v_1(z_1, y)  &\nonumber
   =\Psi(z_1)\Psi(y)^{-1}\bigg[-1-(\la_0+\la_1)\int_{z_1}^y\frac{\Psi(y)\Psi(z_0)^{-1}}{c_0-U'(z_0)}\rmd z_0\bigg]\\
    & =-\beta(z_1, y)-(\la_0+\la_1)\int_{z_1}^y\frac{\beta(z_1, z_0)}{c_0-U'(z_0)}\rmd z_0,\qquad z_0, z_1<y.
\end{align}

Formulae \eqref{eq:m0-sol-1}-\eqref{eq:m1-sol-1} follow from 
 \eqref{eq:mv1} and \eqref{eq:v0-sol}-\eqref{eq:v1-sol}.  


In the case $a_0<x<y=a_+$, by \eqref{eq:bc1} and \eqref{eq:m0m1-1-diff}, 
we obtain the boundary conditions for \eqref{eq:v}: 
\begin{equation*}
v_0(x, y)|_{x\uparrow y}=-1-\la_0/\la_1,\qquad v_1(x, y)|_{x\uparrow y}=0.
\end{equation*}
Similar to \eqref{eq:v0-sol}, we obtain
\begin{equation}
\label{eq:v0-sol=}
v_0(z_0, y)
 =-\beta(z_0, y)\left(1+\la_0/\la_1\right)
 -(\la_0+\la_1)\int_{z_0}^y\frac{\beta(z_0, z_1)}{c_1-U'(z_1)}\rmd z_1,
\end{equation}
and
\begin{equation}
\label{eq:v1-sol=}
v_1(z_1, y)=-(\la_0+\la_1)\int_{z_1}^y\frac{\beta(z_1, z_0)}{c_0-U'(z_0)}\rmd z_0.
\end{equation}
Formulae  \eqref{eq:m0-sol-4} and \eqref{eq:m1-sol-4}
follow from \eqref{eq:mv1} 
and 
\eqref{eq:v0-sol=}-\eqref{eq:v1-sol=}.

The rest of the proof follows by symmetry.
\endproof

\subsection{Large parameters}
If both parameters are large, i.e. $c_0>V$ and $c_1<v$, 
see Fig.\ref{figU'}(d), then 
both potentials, $U(x)-c_0x$ and $U(x)-c_1x$ are single-well,
repulsion points disappear, and the attractors $G_-$ and $G_+$ merge into one $G=[a_-, b_+]$.
In this case, process $X^x(t), x\in G,$ randomly oscillates 
between two attracting stationary states $x=a=a_-$ and $x=b=b_+.$

If both parameters are large, the mean first passage times is derived in closed form using the technique
of Theorem \ref{theo:startingNearAttractor}.

\begin{theorem}
Let the parameters $c_0, c_1$ be large\textup, as described above\textup. 

The average values
$m_0(x,y)=\EE[\sT(x, y)~|~\xi(0)=0]$ and $m_1(x, y)=\EE[\sT(x, y)~|~\xi(0)=1]$ are determined 
as follows\textup:
\begin{itemize}
  \item if $x<y<b$ or $x>y>a,$ then $m_0$ and $m_1$ are determined by formulae 
  \eqref{eq:m0-sol-1}-\eqref{eq:m1-sol-1}\textup;
  \item if $a<x<y<b,$ then formulae \eqref{eq:m0-sol-4}-\eqref{eq:m1-sol-4} are valid\textup;
  \item  if $a<y<x<b,$ then formulae  \eqref{eq:m0-sol-5}-\eqref{eq:m1-sol-6} are valid.
\end{itemize}
\end{theorem}

\section{Stationary measure} \label{sec4}
Let the interval $G=(a, b)$ be invariant under the dynamics $X,$
see \eqref{eq:inv}, so that $c_1\le U'(x)\le c_0,\;x\in G,$
and $U'(b)=c_0,\;U'(a)=c_1$, and thus $G=G_-$ or $G=G_+$, see Fig. \ref{figU'}.

The explicit form of the operator $\sL^*$ adjoint  to $\sL,$ \eqref{def:L},
on the interval $G$ is obtained by integrating by parts
in the integral  
\[
 \langle\sL[\bm f](x),\;\bm g(x)\rangle_{L_2(G)}  
 = \int_G\left[
 (c_0-U'(x))f_0'(x)g_0(x)+ (c_1-U'(x))f_1'(x)g_1(x)\right]\rmd x
\] 
for any test-functions $\bm f=(f_0,\;f_1)$ and $\bm g=(g_0,\;g_1).$
We have
\[
\begin{aligned}
 \langle\sL[\bm f]&(x),\;\bm g(x)\rangle_{L_2(G)}     
 =\left[ (c_0-U'(x))f_0(x)g_0(x)\right]|_{\pd G}      
 +\left[ (c_1-U'(x))f_1(x)g_1(x)\right]|_{\pd G}\\
    -&\int_G\left\{f_0(x)\frac{\rmd}{\rmd x}\left[(c_0-U'(x))g_0(x)\right]
  +f_1(x)\frac{\rmd}{\rmd x}\left[(c_1-U'(x))g_1(x)\right]
    \right\}\rmd x.
\end{aligned}
\]
Therefore, the formal adjoint to the infinitesimal generator $\sL$ is given by
\begin{equation}
\label{eq:L*}
\sL^*\bm g(x)
=\begin{pmatrix}
-  \frac{\rmd}{\rmd x}[(c_0-U'(x))g_0(x)]    &  0  \\ \\
    0  &   - \frac{\rmd}{\rmd x}[(c_1-U'(x))g_1(x)] 
\end{pmatrix}
\end{equation}
with a singular part corresponding to the boundary condition
\begin{equation}
\label{eq:L*sing}
g_0(a)=0,\qquad g_1(b)=0.
\end{equation}

As we show below, there exist invariant measures $\mu^-$ and $\mu^+$ 
 supported on the attractors
 $G_-=(a_-,\;b_-)$ and $G_+=(a_+,\;b_+),$
 respectively. 
\begin{theorem}\label{theo:sm}
Let $G=(a, b)$ be an invariant attractor $(G_-$ or $G_+).$ 
The invariant measure with support $G$ exists and is determined 
by the probability density functions $\pi_0(x)$ and $\pi_1(x):$
 \[
 \pi_0(x)=C_0 \frac{\Psi(x)^{-1}}{c_0-U'(x)},  
 \qquad
   \pi_1(x)  =C_1 \frac{\Psi(x)^{-1}}{U'(x)-c_1},
   \quad x\in G.
\]
Here
$ \Psi(x)^{-1}=\exp\left(-\lambda_0\Phi_0(x)-\lambda_1\Phi_1(x)\right),$ 
and $C_0,\;C_1$ are normalising constants\textup,
\begin{equation}
\label{eq:C0C1}
\aligned
C_0 =\frac{\la_1}{\la_0+\la_1}
\bigg(\int_{a}^{b}\frac{\Psi(x)^{-1}\rmd x}{c_0-U'(x)}\bigg)^{-1},\qquad
C_1=\frac{\la_0}{\la_0+\la_1}
\bigg(\int_{a}^{b}\frac{\Psi(x)^{-1}\rmd x}{U'(x)-c_1}\bigg)^{-1}.
\endaligned
\end{equation}
\end{theorem}
\proof
The invariant density 
$\bm\pi=(\pi_0(x),\;\pi_1(x)),\;x\in G,$ which is defined by
\[
\PP\{X(t)\in\rmd x,\;\xi(t)=0\}\equiv\pi_0(x)\rmd x,\qquad
\PP\{X(t)\in\rmd x,\;\xi(t)=1\}\equiv\pi_1(x)\rmd x,
\]
follows the Fokker-Planck equation 
   \begin{equation*}
(\Lambda^*+\sL^*)\bm\pi=\bm0,\qquad x\in G,
\end{equation*}
$\Lambda^*=\Lambda^\ttT,$ with boundary conditions
$ \pi_0(a)=0$ and $\pi_1(b)=0,$ see  \eqref{eq:L*}-\eqref{eq:L*sing}. 
That is, the entries $\pi_0$ and $\pi_1$ are positive solutions to the coupled equations
\begin{align}
\label{eq:pi0}
-\frac{\rmd}{\rmd x}\big[ (c_0-U'(x))\pi_0(x)\big]
-\la_0\pi_0(x)+\la_1\pi_1(x)=0, &   \\
\label{eq:pi1}
-\frac{\rmd}{\rmd x}\big[ (c_1-U'(x))\pi_1(x)\big]
    +\la_0\pi_0(x)-\la_1\pi_1(x)=0,& 
\end{align}
$ a<x<b.$

 Summing up \eqref{eq:pi0} and \eqref{eq:pi1} we obtain
\begin{equation*}
(c_0-U'(x))\pi_0(x)+(c_1-U'(x))\pi_1(x)\equiv const.
\end{equation*}
The constant is zero.
It follows, for example, from  $U'(a)=c_1,$ and 
due to the boundary conditions, \eqref{eq:L*sing}, $\pi_0(a)=0.$

Therefore,
\begin{equation}
\label{eq:pi0/pi1-symm}
\frac{\pi_0(x)}{\pi_1(x)}=\frac{U'(x)-c_1}{c_0-U'(x)}.
\end{equation}

Equation \eqref{eq:pi0} can be rewritten as
\[
-\frac{[(c_0-U'(x))\pi_0(x)]'}{(c_0-U'(x))\pi_0(x)}-\frac{\la_0}{c_0-U'(x)}
+\frac{\la_1}{c_0-U'(x)}\cdot\frac{\pi_1(x)}{\pi_0(x)}=0.
\]
By \eqref{eq:pi0/pi1-symm}, we have
\[
\frac{\rmd}{\rmd x}\log[(c_0-U'(x))\pi_0(x)]=-\psi(x),\qquad x\in G=(a, b),
\]
where $\psi(x)$ is defined before, \eqref{def:K}, and $\pi_0(a)=0.$ 
 Integrating, we obtain
\begin{equation*}
\pi_0(x)=C_0\frac{\exp(-\la_0\Phi_0(x)-\la_1\Phi_1(x))}{c_0-U'(x)},\qquad a<x<b,
\end{equation*}
where $\Phi_0$ and $\Phi_1$ are defined in \eqref{def:Phii}. The boundary condition
holds, since \[\Phi_0|_{x\to a}<\infty,\qquad \Phi_1(x)|_{x\to a}=+\infty.\]

Similarly,
\begin{equation*}
\pi_1(x)=C_1\frac{\exp(-\la_0\Phi_0(x)-\la_1\Phi_1(x))}{U'(x)-c_1},\qquad a<x<b.
\end{equation*} 

Therefore, $C_0$ and $C_1$ are given by \eqref{eq:C0C1}, 
since the invariant distribution of $\bm\xi(t)$ is determined by
\[
\int_a^b\pi_0(x)\rmd x=\PP_0\{\xi(t)=0\} =\frac{\la_1}{\la_0+\la_1},      \qquad
\int_a^b\pi_1(x)\rmd x=\PP_1\{\xi(t)=0\} =\frac{\la_0}{\la_0+\la_1}. 
\] 
\endproof
\section*{Discussion}
The mathematical interpretation of physical phenomena often uses very sophisticated and 
sometimes  very artificial tools.
For example, a very popular climate change model, widely studied in recent decades,
is based on stochastic differential equation with periodic potential.
\begin{equation}
\label{eq:X*}
\rmd X_\ep(t)=-U'(X_\ep(t), t/T_\ep)\rmd t+\sigma(T_\ep)\rmd W(t),\qquad t>0.
\end{equation}
Here the period $T_\ep$ is usually assumed to be exponentially large in $\ep$
and the periodic fluctuations of the potential are due to periodic changes in the parameters 
of the Earth's orbit
(the so called Milankovitch astronomical cycles, \cite{1930}).

The key property postulated for this model is stochastic resonance.
This methodology involves choosing the diffusion coefficient  $\sigma(T_\ep)$
taking into account the spectral properties of the mean solution of \eqref{eq:X*}.
In my opinion, this method looks artificial, since the using of white noise 
for such modelling (especially for climate models) is not confirmed.
Moreover, it is worth noting that the global and long-term behaviour of the atmosphere 
does not undergo abrupt and rapid changes, and the climate usually changes very slowly.
This makes the classical diffusion process with 
infinite variation and infinitely fast propagation to be problematic for describing 
  climate change.

We propose to add to the (deterministic) continuous periodic changes of potential
randomly atlternating trends corresponding to internal processes 
in the atmosphere and at the surface, see \eqref{eq:X}.
In other words, one can the modify model \eqref{eq:X*}
replacing the white-noise term $\sigma(T_\ep)\rmd W(t)$ by a telegraph process $\TT(t).$

When and if this trend coincides with the direction of the periodic change in potential,
the system can go into a metastable state
receiving a time window for the interwell transition. 
Otherwise, this internal trend has only a retarding effect on the global changes 
that occur from the laws of celestial mechanics.


%
\section*{Appendix: the proof of Theorem \ref{theo:generator}}
Let $\tau^0\sim\mathrm{Exp}(\la_0)$ be the time of the first switching of the 
process $\xi(t)$ from state 0 to state 1,
and $\tau^1\sim\mathrm{Exp}(\la_1)$ be the time of the first switching from state 1 to state 0.

Let 
$\cP(t, \rmd y~|~x)=\left(p_{ij}(t, \rmd y~|~x)\right)_{i, j\in\{0, 1\}}$  
be the transition probability matrix,
where \[p_{ij}(t, \rmd y~|~x)=\PP\{X(t)\in\rmd y, \xi(t)=j~|~X(0)=x, \xi(0)=i\}.\]
The corresponding Markov semigroup $P_t$ 
is defined on the test-function $\bm f=(f_0, f_1)$ as follows:
\[
(P_t\bm f)_i(x)=\EE\left[\bm f(\Xi(t))~|~\Xi(0)=(x, i)\right]
=\int_{-\infty}^\infty f_0(y)p_{i0}(t, \rmd y~|~x)
+\int_{-\infty}^\infty f_1(y)p_{i1}(t, \rmd y~|~x),
\]
$i\in\{0, 1\}.$

Conditioning on the first switching, see \eqref{eq:conditioning}, we obtain the following system 
of integral equations:  
\begin{equation}
\label{eq:int}
\left\{
\aligned
p_{0j}(t, \rmd y~|~x)=&\rme^{-\la_0t}\delta_{\gamma_0(t, x)}(\rmd y)
+\int_0^t\la_0\rme^{-\la_0\tau}p_{1j}(t-\tau, \rmd y~|~\gamma_0(\tau, x))\rmd \tau,\\
p_{1j}(t, \rmd y~|~x)=&\rme^{-\la_1t}\delta_{\gamma_1(t, x)}(\rmd y)
+\int_0^t\la_0\rme^{-\la_0\tau}p_{0j}(t-\tau, \rmd y~|~\gamma_1(\tau, x))\rmd \tau,\\
&\qquad j\in\{0, 1\}.
\endaligned
\right.
\end{equation}
Here $\delta_{a}(\rmd y)$  is Dirac's $\delta$-measure at point $a$.

By applying the operator $ \dfrac{\pd}{\pd t}-L_0^x  $ to the first equation of \eqref{eq:int}, 
and the operator
$\dfrac{\pd}{\pd t}-L_1^x$ to the second,
one can find that system of integral equations \eqref{eq:int}
is equivalent to some initial value problem for partial differential equations.  
Indeed, by virtue of \eqref{eq:Lgamma},
\[
\left\{
\begin{aligned}
  \frac{\pd p_{0j}(t ,\rmd y~|~x)}{\pd t}&-L_0^x[p_{0j}(t ,\rmd y~|~x)]  
   =-\la_0\rme^{-\la_0t}\delta_{\gamma_0(t, x)}(\rmd y) 
    \\
    + \la_0&\rme^{-\la_0t}p_{1j}(0, \rmd y~|~\gamma_0(t, x))
      -\int_0^t\la_0\rme^{-\la_0\tau}\frac{\pd}{\pd \tau}
   \big[ p_{1j}(t-\tau, \rmd y~|~\gamma_0(\tau, x))\big]\rmd \tau,\\
  \frac{\pd p_{1j}(t ,\rmd y~|~x)}{\pd t}&-L_1^x[p_{0j}(t ,\rmd y~|~x)]  
   =-\la_1\rme^{-\la_1t}\delta_{\gamma_1(t, x)}(\rmd y) 
    \\
     + \la_1&\rme^{-\la_1t}p_{0j}(0, \rmd y~|~\gamma_1(t, x)) -\int_0^t\la_1\rme^{-\la_1\tau}
    \frac{\pd}{\pd \tau}  \big[  p_{0j}(t-\tau, \rmd y~|~\gamma_1(\tau, x))\big]\rmd \tau.
\end{aligned}
\right.\]
Integrating by parts, we obtain
\begin{equation*}
\left\{
\begin{aligned}
   \frac{\pd p_{0j}(t ,\rmd y~|~x)}{\pd t}&-L_0^x[p_{0j}(t ,\rmd y~|~x)]  
   =       -\la_0p_{0j}(t ,\rmd y~|~x)+\la_0p_{1j}(t ,\rmd y~|~x),\\
   \frac{\pd p_{1j}(t ,\rmd y~|~x)}{\pd t}&-L_1^x[p_{1j}(t ,\rmd y~|~x)]  
   =     \la_1p_{0j}(t ,\rmd y~|~x)  -\la_1p_{1j}(t ,\rmd y~|~x),
   \\
&\qquad j\in\{0, 1\},
\end{aligned}
\right.\end{equation*}
which can be rewritten in matrix form
\begin{equation*}
\frac{\pd \cP(t, \rmd y~|~x)}{\pd t}=\left(\Lambda+\sL^x\right)[\cP(t, \rmd y~|~x)],\qquad t>0,
\end{equation*}
with initial condition
\[
\left.\cP(t, \rmd y~|~x)\right|_{t\downarrow0}=
\begin{pmatrix}
    \delta_x(\rmd y)  &  0  \\
   0   &  \delta_x(\rmd y)
\end{pmatrix}.
\]\endproof

The research was supported by the Russian Science Foundation (RSF), 
project number 24-21-00245,  https://rscf.ru/project/24-21-00245

\section*{Acknowlegements}
This research was inspired by a report given at the 
United Seminar of the Department of Probability Theory of Moscow State University,
by Professor A.N.Shiryaev, \\
https://www.youtube.com/watch?v=KWzp5ruOSP4

\noindent 
I am very grateful for opportunity to participate in this seminar.

I am also grateful to the reviewers and the Associate Editor 
for their careful reading of the manuscript and 
very helpful suggestions that made the text much better.

\end{document}